\documentclass[11pt,twoside,reqno]{amsart}
\usepackage[all]{xy}
        \usepackage {amssymb,latexsym,amsthm,amsmath,mathtools}
        \usepackage{enumitem,color}
        
\usepackage{url}
\usepackage[colorlinks,citecolor=red,urlcolor=blue,bookmarks=false,hypertexnames=true]{hyperref}
\usepackage{mathtools}

\usepackage[backend=biber]{biblatex}
\long\def\symbolfootnote[#1]#2{\begingroup%
\def\thefootnote{\fnsymbol{footnote}}\footnote[#1]{#2}\endgroup}

\makeatletter
\def\imod#1{\allowbreak\mkern10mu({\operator@font mod}\,\,#1)}
\makeatother

\newtheorem{theorem}{Theorem}[section]
\newtheorem{lemma}[theorem]{Lemma}

\newtheorem*{theorem*}{Theorem}
\theoremstyle{definition}
\newtheorem{definition}[theorem]{Definition}

\newtheorem{example}[theorem]{Example}
\numberwithin{equation}{section}

\newcommand{\ignore}[1]{}

\newcommand{\mynote}[1]{}

\newcommand{\Irr}{\textnormal{Irr}}
\newcommand{\cd}{\textnormal{cd}}

\newcommand{\Core}{\textnormal{Core}}

\newcommand{\Gal}{\textnormal{Gal}}

\newcommand{\rdim}{\textnormal{rdim}}
\newcommand{\ed}{\textnormal{ed}}
\newcommand{\irr}{\textnormal{irr}}

\begin{document}
\setcounter{section}{0}
\title{Various Representation Dimensions associated with a Finite Group}
\author{Anupam Singh}
\address{IISER Pune, Dr. Homi Bhabha Road, Pashan, Pune 411 008, India}
\email{anupamk18@gmail.com}
\author{Ayush Udeep}
\address{IISER Mohali, Knowledge city, Sector 81, Mohali 140 306, India}
\email{udeepayush@gmail.com}
\thanks{The first-named author is funded by an NBHM research grant 
02011/23/2023/NBHM(RP)/RDII/ 5955 for this research. The second-named author thanks IISER Mohali, India, for a postdoctoral fellowship.}


\begin{abstract}
To a finite group $G$, one can associate several notions of dimensions (or degrees). In this survey, we attempt to bring together some of the notions of dimensions or degrees defined using representations of the group in General Linear Groups and permutation groups. These are embedding degree, minimal faithful irreducible character degree, minimal faithful permutation representation degree, minimal faithful quasi-permutation representation degree and essential dimension. We briefly present the progress in understanding these notions and the related problems. 
\end{abstract}

\maketitle

\section{Introduction}
Let $G$ be a finite group. As soon as we are introduced to a group via axiomatic definition in our undergraduate courses, several questions come to mind. One among them is can we think of a finite group as a subgroup of some concrete group? Indeed an answer is provided by Cayley's theorem that a finite group $G$ of size $n$ can be embedded inside $S_n$ using the left multiplication map. We can also use a similar method to embed $G$ inside $GL_n(\mathbb C)$ giving rise to the regular representation. The underlying concept is that we supposedly understand the structure of symmetric groups and linear groups better which involve combinatorics and linear algebra, respectively. One can further ask questions the other way around that for a fixed $d$ and a field $k$ what are the finite subgroups of a given $GL_d(k)$ and a $S_d$? These questions date back to Minkowaski (1887) and  Jordan (1878) which is beautifully described by Serre in his book Finite Groups: an Introduction~\cite{S22}.  

Indeed, these ideas led to the definition of several different kinds of ``dimensions'' or ``degrees'' associated with a group $G$. In this note, we briefly explore what is embedding degree $\delta(G)$, minimal faithful irreducible character degree $\delta_{irr}(G)$, minimal faithful permutation representation degree $\mu(G)$, minimal faithful quasi-permutation representation degree $q(G)$, and essential dimension $ed_K(G)$. These notions are studied in literature some of which we aim to bring together. We explain the concepts through some easy examples and provide updates on the related problems. Throughout this article, $p$ denotes a prime, and $\mathbb{C}$ is the underlying field unless stated otherwise.

\subsection*{Acknowledgment} This note is based on a talk given by the first named author at Harish-Chandra Research Institute (HRI), Prayagraj in a workshop organised by Professor Manoj Kumar where the second author was also a participant. Both authors express their gratitude to the institute for its wonderful hospitality.  

\section{Embedding Degree} \label{sec:m(G)}
The \emph{embedding degree}, or the representation dimension, of a finite group $G$, denoted $\delta(G)$ or by $\rdim(G)$, is the minimal dimension of a faithful complex representation of $G$. Alternatively, $\delta(G)$ is the smallest positive integer $n$ such that $G$ embeds into $GL_n(\mathbb{C})$, i.e., $G$ is isomorphic to a linear group of degree $n$. The determination of the embedding degree of a group is one of the classical and challenging problems in the representation theory of finite groups, and its study has found many applications, for example, its connections with the \emph{essential dimension}, $\ed(G)$, of a group $G$. We discuss the essential dimension of $G$ in Section~\ref{sec:ed(G)}. Note that the determination of all finite groups with representation dimension $\delta$ is equivalent to the determination of all the finite subgroups of $GL_{\delta}( \mathbb{C})$.

Motivated by the relation between the embedding degree and the essential dimension, in 2011, Cernele et al.~\cite{CKR11} initiated the study of the embedding degree and obtained an upper bound for $\delta(G)$ of a $p$-group $G$, where $p$ is a prime. They prove that if the order of $G$ is $p^n$ and the rank of $\mathcal Z(G)$ is $r$, then for almost all pairs $(p, n)$, the maximum value of $\delta(G)$ is $f_{p}(n)$, where
$$f_{p}(n) = \max_{r \in \mathbb{N}}\left( rp^{\lfloor (n-r)/2 \rfloor} \right).$$
The exceptional cases are of $(p, n)$ being $(2, 5)$, $(2, 7)$ and $(p, 4)$ where $p$ is odd.
In these cases the maximum value of $\delta(G)$ is $5, 10$, and $p+1$, respectively. Note that,
\begin{example}
If $G$ is a cyclic group then $\delta(G) = 1$. In fact, $\delta(G)=1$ if and only if $G$ is cyclic. However, if $G \cong C_2 \times C_2$, then $\delta(G) = 2$. 
\end{example}
\begin{example}
When $G$ is dihedral group $D_n$, $Q_8$, and $C_n \times C_m$ where $n$ is not coprime to $m$ then $\delta(G)=2$. Note that $\delta(G)=2$ amounts to finding all finite subgroups of $GL_2(\mathbb C)$. 
\end{example}
\noindent This shows that the knowledge of the character table is not sufficient for the computation of the embedding degree of a group. Utilizing a result of \cite{MR10}, Bardestani et al. \cite{BMS16} proved that if $G$ is a $p$-group and $\rho$ is a faithful representation of $G$ of least dimension, then $\rho$ decomposes into exactly $r$ irreducible representations $\rho_i ~(1\leq i \leq r)$ of $G$, where $r$ is the rank of $\mathcal Z(G)$. Then
$$ 1 = \ker(\rho) = \bigcap_{i=1}^{r} \ker(\rho_{i}) \, \text{ and } \bigcap_{i=1, i\neq j}^{r} \ker(\rho_{i}) \neq 1 \text{ for each } j ~(1\leq j \leq r).$$
Thus, $\delta(G)$ may be restated as
\begin{align*}
\delta(G) = \min \bigg\{  & \sum_{\chi \in X} \chi(1) : X \subset \Irr(G) \text{ such that } \bigcap_{\chi \in X} \ker(\chi) = 1  \\
 	&  \text{ and } \bigcap_{\chi \in Y} \ker(\chi) \neq 1 \text{ for all } Y \subsetneq X \bigg\}.
 \end{align*}
Bardestani et al. \cite{BMS16} further computed the embedding degree of the Heisenberg group over $\mathcal{O}/\mathcal{I}$, where $\mathcal{O}$ is a certain ring of integers and $\mathcal{I}$ is its maximal ideal. Moret\'{o} \cite{M21} proved that for a finite group $G$, $\delta(G) \leq \sqrt{|G|}$ or $\frac{3}{\sqrt{8}}\sqrt{|G|}$. Further, Moret\'{o} classified groups for which $\delta(G) = \sqrt{|G|}$; these groups turn out to be $2$-groups. It also proved that the embedding degree of a finite abelian group is the rank of the group. This implies that the embedding degree for finite groups may be arbitrarily large. In \cite{KKS23}, Kaur et al. computed the embedding degree for various classes of groups, such as some finite $p$-groups ($p$ is a prime), certain direct products of finite groups, odd order groups with exactly two nonlinear irreducible characters of each degree, and finite groups whose all nonlinear irreducible characters are of distinct degree. They also provide {\sc GAP} \cite{GAP} codes for the computation of the embedding degree. 
L\"{u}beck \cite{L01} and Tiep and Zalesskii \cite{TZ96} have studied representation dimension for finite groups of Lie type.

\section{Minimal Faithful Irreducible Character Degree} \label{sec:mirr(G)}

For a group $G$, the \emph{minimal faithful irreducible character degree} of $G$, $\delta_{\irr}(G)$, is defined as
\[ \delta_{\irr}(G) := \min \{ \rho(1) : \rho \text{ is an irreducible faithful representation of } G \}. \]
Note that a group may not have a faithful irreducible character and hence $\delta_{\irr}(G)$ may not exist. 
\begin{example}
For $G=C_2\times C_2$, $\delta_{irr}(G)$ does not exist. 
\end{example}
\noindent The study of the faithful irreducible representations of a group has been an intriguing area of research. Since the beginning of $20^{\text{th}}$ century, many researchers have determined various characterizations for groups which possess a faithful irreducible character; Szechtman provides a detailed history of the problem in~\cite{S16}. Kitture and Pradhan \cite{KP22} obtained the degrees of faithful irreducible representations for a class of metabelian groups. When a faithful irreducible character of $G$ exists (say $\chi$), we have
the inequality $\delta(G) \leq \chi(1)$. However, the quantity $\delta(G)$ may be different from $\delta_{\irr}(G)$. 
\begin{example} For the group  $A_4 \times D_{10}$, we have $\delta(A_4 \times D_{10}) = 5$ and $\delta_{\irr}(A_4 \times D_{10}) = 6$ (see Example 2.4 of \cite{KKS23}).
\end{example}

In the literature, the minimal faithful irreducible character degree for various classes of groups have been studied. In the case of non-abelian finite simple groups, every nontrivial character is faithful. Hence, for these groups, the embedding degree is equal to the minimal degree of
a nonlinear irreducible character. By Lemma 3.1 of \cite{BTVZ17}, we get $\delta(A_{n}) = \delta_{\irr}(A_{n}) = n-1$ for $n\geq 15$. It is easy to compute from {\sc GAP} that $\delta(A_n) = n-1$ for $6\leq n \leq 14$, and $\delta(A_5) = 3$. Further, from a result of Rasala \cite{R77}, we get $\delta_{\irr}(S_n) = n-1$ for $n\geq 5$. Since the standard character of $S_{n}$ is faithful and has degree $n-1$, it follows that $\delta(S_{n}) = n-1$ for $n\geq 5$. From Steinberg's work in \cite{S51}, we get that $\delta_{\irr}(GL(2, q)) = q-1$, $\delta_{\irr}(GL(3, q)) = q^2+q$ and $\delta_{\irr}(GL(4, q)) = (q +1)(q^2 +1)$. In \cite{M99}, Mann proved that all faithful irreducible characters of a normally monomial group are of the same degree, which is the maximum degree of all irreducible characters of the group. As a result, if $G$ is a normally monomial $p$-group with cyclic center, then $\delta_{\irr}(G) = \max \cd(G)$, where $\cd(G)$ denotes the character degree set of $G$. Kaur et al. \cite{KKS23} proved that if $G$ is a finite nilpotent group with the cyclic centre, then $\delta(G) = \delta_{\irr}(G)$. Hence, we also obtain the embedding degree for any normally monomial $p$-group with a cyclic centre from Mann's result. In \cite{KKS23}, Kaur et al. obtained a condition when $\delta(G_1 \times G_2) = \delta_{\irr}(G_1 \times G_2) = \delta_{\irr}(G_1) \delta_{\irr}(G_2)$, where $G_1$ and $G_2$ are finite non-abelian groups.\\

\noindent {\bf Question 1.} Characterize groups $G_1$ and $G_2$ such that $\delta(G_1 \times G_2) = \delta_{\irr}(G_1 \times G_2) $.\\

Kaur et al. also computed $\delta(G)$ and $\delta_{\irr}(G)$ where $G$ is a finite non-abelian group of odd order with exactly two nonlinear irreducible characters of each degree or a finite non-abelian group whose all nonlinear irreducible characters have distinct degrees.

\section{Minimal Faithful Permutation Representation Degree} \label{sec:mu(G)}

Cayley's theorem states that every group is isomorphic to a subgroup of a symmetric group. It is interesting to find the least positive integer $n$ such that $G$ is embedded in the symmetric group $S_{n}$. Such $n$ is called the \emph{minimal faithful permutation representation degree} of $G$, denoted by $\mu(G)$. One can easily observe that for a group $G$, we have $\delta(G) \leq \mu(G)$. It is well known that the minimal faithful permutation representation degree of $G$ is given by
\[  \mu(G)= \min \left\{ \sum_{i=1}^{k}|G:G_{i}| : G_{i}\leq G \text{ for } 1\leq i \leq k \text{ and } \bigcap_{\substack{i=1}}^{k} \Core_{G}(G_{i}) =1  \right\}, \]
where $\Core_{G}(G_i)$ is the core of $G_i$ in $G$.
Here $\mathcal{H}=\{G_{1},\ldots,G_{k}\}$ is called a \emph{minimal faithful permutation representation} of $G$. Let $H$ be a group and $\{ H_1, H_2, \ldots, H_s \}$ be a minimal faithful permutation representation of $H$. Then one can show that
\[ \{ G_1 \times H, \ldots, G_r \times H, G \times H_1, \ldots, G \times H_s \} \]
is a minimal faithful permutation representation of $G\times H$. \\
\begin{example}
One can show the following: $\mu(C_p)=p$, $p$ a prime; $\mu(C_3\times C_3) = 6 = \mu(C_3) + \mu(C_3)$; $\mu(C_6=5)$; $\mu(D_4)=4$ ($|D_4|=8$); $\mu(Q_8)=8$.  
\end{example}

\noindent Karpilovsky \cite{K70} first computed $\mu(G)$ for a finite abelian group $G$. Johnson \cite{J71} proved that $\mu(G) = |G|$ if and only if $G$ is isomorphic to a cyclic $p$-group (for a prime $p$), a generalized quaternion group or the non-cyclic group of order
$4$. Johnson also computed the cardinality of a minimal faithful permutation representation of a $p$-group. Various researchers investigated the relation between $\mu(G\times H)$ and $\mu(G)+\mu(H)$ for two groups $G$ and $H$. Johnson \cite[Proposition 2]{J71} proved that for any two groups $G$ and $H$, \[ \mu(G\times H)\leq \mu(G)+ \mu(H), \]
and the reverse inequality holds whenever $G$ and $H$ have coprime orders. Wright \cite{W75} improved this result and proved that for all non-trivial $G, H \in$ \{nilpotent, symmetric, alternating or dihedral groups\}, $\mu(G\times H)= \mu(G)+\mu(H)$. Wright further provided an example of groups $G$ and $H$ such that $\mu(G\times H)\neq \mu(G)+\mu(H)$ (see \cite[Section 5]{W75}). Easdown and Praeger \cite{EP88} proved that for finite simple groups $G$ and $H$, $\mu(G\times H)= \mu(G)+\mu(H)$.
Saunders \cite{S09} produced an infinite family of examples of permutation groups $G$ and $H$ where $\mu(G\times H) < \mu(G) + \mu(H)$. In \cite{S10}, Saunders computed the minimal faithful permutation degree of a class of finite reflection groups and proved that they form examples where the minimal degree of a direct product is strictly less than the sum of the minimal degrees of the direct factors. Easdown and Saunders \cite{ES16} further proved that $10$ is minimal in the sense that $\mu(G\times H) = \mu(G) + \mu(H)$ for all groups $G$ and $H$ such that $\mu(G\times H) \leq 9$. Easdown and Hendriksen \cite{EH16} proved that if $G$ and $H$ are finite groups, then $\mu(G\rtimes H) \leq |G| + \mu(H)$.

Several researchers, for example, Babai et al. \cite{BGP93} and Franchi \cite{F11} investigated bounds of $\mu(G)$ for a group $G$. Babai et al. \cite{BGP93} proved that $\nu(G) \geq |G|/\mu(G) \geq f(\nu(G))$, where $\nu(G)$ is the index of the largest cyclic subgroup of prime-power order in $G$ and $f(n) = e^{(c\sqrt{\log n})}$ for a constant $c>0$. Franchi \cite{F11} proved that if $G$ is a finite $p$-group with an abelian maximal subgroup, then $\mu(G/G') \leq \mu(G)$. Berkovich \cite{YB99} has obtained various results by investigating the relation between $\mu(G)$ and $i(G)$ for a group $G$, where $i(G) = \min\{ |G:H| ~:~ H\leq G \}$ if $G> 1$ and $i(G) = 1$ if $G=1$. Hall and Senior \cite{HS64} listed minimal degree faithful representations for the groups of order dividing $2^6$.
Lemieux \cite{LThesis} determined the values of $\mu(G)$ for the groups of order dividing $p^4$.

For a group $G$ and its normal subgroup $N$, various researchers have investigated the relation between $\mu(G)$ and $\mu(G/N)$. In \cite{N86}, Neumann took $G$ to be the direct product of $k>1$ copies of the dihedral group of order $8$ to show that $\mu(G/N) > \mu(G)$ for some normal subgroup $N$ of $G$. Easdown and Praeger \cite{EP88} produced further examples of such groups $G$ with quotient groups $G/N$ such that $\mu(G/N) > \mu(G)$. A group $G$ is called exceptional if there exists a normal subgroup $N$ such that $\mu(G/N) > \mu(G)$; here $N$ is called a distinguished subgroup and $G/N$ is called a distinguished quotient. Kov\'{a}cs and Praeger proved that $\mu(G/N) \leq \mu(G)$ if $G/N$ has no nontrivial abelian normal subgroup. Holt and Walton \cite{HW02} proved that there exists a constant $c$ such that $\mu(G/N) \leq c^{\mu(G) - 1}$ for all finite groups $G$ and all normal subgroups $N$ of $G$. In \cite{LThesis, SL07}, Lemieux proved that there are no exceptional $p$-groups of order less than $p^5$, and also produced an exceptional group of order $p^5$ where $p$ is an odd prime. Chamberlain \cite{C18} independently proved that there are no exceptional p-groups of order less than $p^5$. In \cite{S14}, Saunders computed the minimal faithful permutation degree of the irreducible Coxeter groups and exhibited examples of new exceptional groups. Britnell et al. \cite{BSS17} classified all exceptional $p$-groups of order $p^5$ where $p$ is a prime and proved that the proportion of groups of order $p^5$ that are exceptional is asymptotically $1/2$. Due to this observation, they emphasize that the established term exceptional for these groups may be less appropriate. However, they have not commented on the case of groups of order $p^6$.

\section{Minimal Faithful Quasi-permutation Representation Degree} \label{sec:c(G)}
A quasi-permutation matrix is defined to be a complex square matrix with a non-negative integral trace. In 1963, Wong~\cite{W63} defined a Quasi-permutation group to be a finite group in which every element is a quasi-permutation matrix. The above terminology allows us to have two more degrees. 
\begin{definition}
The minimal degree of a faithful representation of $G$ by quasi-permutation matrices over the field $\mathbb{Q}$ is denoted by $q(G)$. The minimal degree of a faithful representation of $G$ by complex quasi-permutation matrices is denoted by $c(G)$.
\end{definition}
\noindent Since every permutation matrix is a quasi-permutation matrix, it is easy to see that $c(G)\leq q(G)\leq \mu(G)$, i.e., $c(G)$ and $q(G)$ provide a lower bound for $\mu(G)$. The above inequalities could be strict, as well as, equality. 
\begin{example}
Take $G_{1}=SL(2,5)$, $G_{2} = C_{6}$, $G_{3} = Q_{8}$ (the quaternion group of order 8), and $G_{4} = D_{4}$ (the dihedral group of order $8$). Then $c(G_{1})=8$, $q(G_1)=16$ and $\mu(G_1)=24$ (see \cite{BGHS94}), $c(G_2) = q(C_2) = 4$ and $c(G_2) = 5$, $c(G_3) = 4$ and $q(G_3) = \mu(G_3) = 8$, and $c(G_4) = q(G_4) = \mu(G_4) = 4$. Of course, we have $\delta(G) \leq c(G)$.
\end{example}

\noindent Burns et al. \cite{BGHS94} initiated the computation of $c(G)$ and $q(G)$ for abelian groups and proved that if $G$ is a finite abelian group, then $c(G) = q(G) =\mu(G)$ if and only if $G$ has no direct factor of order $6$. Behravesh \cite{B97b} improved the above result and proved that if
$G \cong \prod_{i = 1}^k C_{p^{a_{i}}}$ and $n$ is maximal such that $G$ has a direct factor isomorphic to $C_{6}^{n}$, then $c(G) = q(G)=\mu(G)-n$. Behravesh \cite{B97a} provided an algorithm for computing $c(G)$ and $q(G)$ for a finite group $G$. Behravesh and Ghaffarzadeh \cite{BG11} improved the said algorithm and proved the following: 
\begin{lemma} \textnormal{\cite[Lemma 2.2]{BG11}}\label{lemma:c(G)Algorithm}
Let $G$ be a finite group.  Let $X \subset \Irr(G)$ be such that $\bigcap_{\chi \in X} \ker (\chi)= 1$ and $\bigcap_{\chi \in Y} \ker (\chi) \neq 1$ for every proper subset $Y$ of $X$. Define $\xi_X = \sum_{\chi \in X} \sum_{\sigma \in \Gamma(\chi)} \chi^{\sigma}$ and $\xi'_X = \sum_{\chi \in X} \left[ m_{\mathbb{Q}}(\chi) \sum_{\sigma \in \Gamma(\chi)} 	\chi^{\sigma}\right]$, where for $\chi\in \Irr(G)$, $\Gamma(\chi) = \Gal(\mathbb{Q}(\chi): \mathbb{Q})$ and $m_{\mathbb{Q}}(\chi)$ is the Schur index of $\chi$ over $\mathbb{Q}$. Let $m(\xi_X)$ and $m(\xi'_X)$ be the absolute value of the minimum value that $\xi_X$ and $\xi'_X$ take over $G$ respectively. Then
\begin{align*}
c(G) =& \min \{\xi_X(1) + m(\xi_X) \; | \; X \subset \Irr(G)  \text{\ satisfying\ the\ above\ property} \}, \text{ and} \\
q(G) =& \min \{\xi'_X(1) + m(\xi'_X) \; | \; X \subset \Irr(G)  \text{\ satisfying\ the\ above\ property} \}.
	\end{align*}  
\end{lemma}
\noindent The Proof of Lemma \ref{lemma:c(G)Algorithm} contains a minor error, and Prajapati and Udeep have included a correct proof in \cite{PU23a}. Behravesh and Ghaffarzadeh \cite{BG11} proved that if $G$ is a $p$-group ($p$ is a prime), then $q(G)=\mu(G)$. Moreover, if $p\neq 2$, then $c(G)= q(G)= \mu(G)$. Hence if $G = H\times K$, where $H$ and $K$ are $p$-groups ($p$ an odd prime), then $c(G) = c(H) + c(K)$ and $q(G) = q(H) + q(K)$. Ghaffarzadeh and Abbaspour \cite{GA12} further proved that the above result is true even for $2$-groups $H$ and $K$. In the same article, they proved that if $X\subset \Irr(G)$ satisfying the hypothesis of Lemma \ref{lemma:c(G)Algorithm} such that $c(G) = \xi_X(1) + m(\xi_X)$, then the rank of $\mathcal{Z}(G)$ is the cardinality of $X$. 

Various researchers have worked out $c(G)$ and $q(G)$ for different classes of groups. In \cite{B97a}, Behravesh proved that if $G$ is a $p$-group of nilpotency class $2$ with cyclic center, then $c(G) = |G/\mathcal{Z}(G)|^{1/2} |\mathcal{Z}(G)|$. Behravesh computed $c(G)$, $q(G)$ and $\mu(G)$ when $G$ is $SL(2,q)$ or $PSL(2,q)$ in \cite{B99}; Darafsheh et al. \cite{DGDB01} computed the degrees for $G = GL(2, q)$, and Darafsheh and Ghorbany \cite{DG03} computed the degrees when $G$ is $SU(3, q^2)$ or $PSU(3, q^2)$.

One can see significant progress in the case of $p$-groups, where $p$ is a prime. Since for $p$-groups $G$ of odd order, $\mu(G) = c(G) = q(G)$, the quantities $c(G)$ and $q(G)$ for groups $G$ of order $p^4$ ($p$ odd) are derived from Lemieux's work \cite{LThesis}. In \cite{B00}, Behravesh computed $c(G)$ when $G$ is a metacyclic $p$-group with a non-cyclic center. Behravesh and Delfani \cite{BD18} claim to list $c(G)$ for the groups of order $p^5$ ($p$ odd); however, their work has many inaccuracies. Prajapati and Udeep have computed $\mu(G)$ and $c(G)$ for groups $G$ of odd order $p^5$ in \cite{PU24}. Pradhan and Sury \cite{PS22} studied the quasi-permutation representations as well as the degrees $c(G)$, $q(G)$ and $\mu(G)$ for the holomorph of cyclic $p$-groups, where $p$ is a prime. In \cite{PU23a, PU23b}, Prajapati and Udeep studied $c(G)$, $q(G)$ and $\mu(G)$ for various classes of non-abelian $p$-groups such as VZ $p$-groups, Camina $p$-groups, and the groups with cyclic center. In \cite{PU23b}, authors prove that if $G$ is either a VZ $p$-group with $d(G') = d(\mathcal{Z}(G))$, or a Camina $p$-group, then $c(G) = |G/\mathcal{Z}(G)|^{1/2}c(\mathcal{Z}(G))$. In \cite{O'BPU23}, O'Brien et al. computed $c(G)$ and $\mu(G)$ for all the groups of order $p^6$ ($p$ odd); they also provide a minimal faithful permutation representation for each group. In \cite{PU23c}, Prajapati and Udeep establish equality between $c(G)$ and a $\mathbb{Q}_{\geq 0}$-sum of codegrees of some irreducible characters of a non-abelian $p$-group $G$ of odd order. Note that for $\chi \in \Irr(G)$, codegree of $\chi$ is defined as $\frac{|G/\ker(\chi)|}{\chi(1)}$. In the same article, Prajapati and Udeep prove that if $G$ is a GVZ $p$-group with a cyclic center, then $c(G) = |G/\mathcal{Z}(G)|^{1/2}|\mathcal{Z}(G)|$.\\

\noindent {\bf Question 2.} Characterize $p$-groups $G$ for which $c(G) = |G/\mathcal{Z}(G)|^{1/2}|\mathcal{Z}(G)|$.

\section{Essential Dimension} \label{sec:ed(G)}
Buhler and Reichstein \cite{BR97} established the notion of essential dimension in the following way. Let $K$ be a field and $L/L_{0}$ be a finite separable extension containing $K$. We say $L/L_{0}$ is defined over $E_{0}$ if there is a field $E$ such that $K \subset E_{0} \subset E \subset L$, $[E:E_{0}] = [L:L_{0}]$ and $L = EL_{0}$.
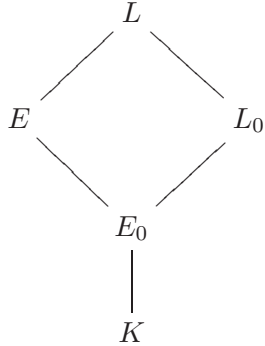
\[
\xymatrix{ &L\ar@{-}[ld]\ar@{-}[rd] \\
	E\ar@{-}[rd] &&L_0 \ar@{-}[ld]\\
	& E_0& \\
	&K\ar@{-}[u] }
\]
The essential dimension of $L/L_{0}$ is given by
\[ ed_K(L/L_0) = \min tr.deg_K E_0 \]
where $E_0$ runs through all the intermediate fields over which $L/L_0$ is defined, and $tr.deg_K$ denotes the transcendence degree of $E_0$ over $K$.
Further, Buhler and Reichstein defined the essential dimension of a finite group, denoted $ed_K(G)$, as the essential dimension of any Noether extension $K(V)/K(V)^G$, where $G\rightarrow GL(V)$ is a faithful representation of $G$ and $K(V)$ is the function field of the affine space $V$ over $K$. The integer $ed_K(G)$ is the smallest number of algebraically independent parameters required to define Galois $G$-algebras over any field extension of $K$. There is another functorial definition of essential dimension due to Merkurjev \cite{M00} in his unpublished notes, which Berhuy and Favi later mentioned in their paper (see \cite{BF03}) too. 

Various researchers have studied essential dimensions in the recent past. Buhler and Reichstein \cite{BR97} proved that $ed_K(H_1 \times H_2) \leq ed_{K}(H_1) + ed_{K}(H_2)$ for groups $H_1$ and $H_2$. They also proved that if $K$ is a field containing all roots of unity, then $ed_K(G) = 1$ if and only if $G$ is isomorphic to a cyclic group or a dihedral group $D_{2m}$ where $m$ is odd. Ledet \cite{L07} further proved that a finite group $G$ has essential dimension $1$ over an infinite field $K$ if and only if there exists an embedding of $G$ in $GL_2(K)$ such that the image of $G$ contains no scalar matrices other than the identity. In \cite{R00}, Reichstein computed essential dimensions for various classes of groups such as Symmetric groups, Orthogonal groups, General linear groups, etc. In \cite{BF03}, Berhuy and Favi proved that for a finite group $G$, $ed_K(G) \leq \delta_{K}(G)$, where $\delta_{K}(G)$ is the embedding degree of $G$ over $K$. In fact, Karpenko and Merkurjev \cite{KM08} proved that if $G$ is a $p$-group and $K$ is a field of characteristic different from $p$ containing a primitive $p^{th}$ root of unity, then $ed_K(G) = \delta_{K}(G)$. Florence \cite{F08} computed the essential dimension of a cyclic $p$-group where $K$ in a field of characteristic not $p$, containing the $p^{th}$ roots of unity. Further, Wong \cite{W11} computed an upper bound of the essential dimension of the cyclic group of order $p_{1}^{n_1}\cdots p_{r}^{n_r}$ over a field $K$ of characteristic different from $p_{i}$ containing all the primitive $p_{i}^{th}$ roots of unity, where $p_{i}$'s are distinct prime numbers. Brosnan et al. \cite{BRV10} studied the essential dimension for the spinor group. Readers can see detailed survey articles by Reichstein \cite{R10} and Merkurjev \cite{M13,M17} on the essential dimension.


\end{document}